\documentclass[a4paper,12pt]{article}
\usepackage{comment}
\usepackage{cite}
\usepackage{amsmath}
\usepackage{amssymb}
\usepackage{amsfonts}
\usepackage[T1]{fontenc}
\usepackage[utf8]{inputenc}
\usepackage{graphicx}
\usepackage{fancyhdr}
\usepackage{float}
\usepackage{xcolor}
\usepackage{authblk}
\usepackage{mathrsfs}
\usepackage{empheq}
\usepackage{url}

\usepackage{geometry}
\begin{document}

\title{Series representations for $\pi^3$ involving the golden ratio}

\author[$\dagger$]{Jean-Christophe {\sc Pain}\\
\small
CEA, DAM, DIF, F-91297 Arpajon, France\\
Universit\'e Paris-Saclay, CEA, Laboratoire Mati\`ere en Conditions Extr\^emes,\\ 
91680 Bruy\`eres-le-Ch\^atel, France
}

\maketitle

\begin{abstract}
Although many series exist for $\pi$ and $\pi^2$, very few are known for $\pi^3$. In this article, we derive, using a trigonometric identity obtained by Euler, two representations of $\pi^3$ involving infinite sums and the golden ratio. The methodology can be generalized in order to obtain further series, relating by the way $\pi^3$ to other mathematical constants.
\end{abstract}

\section{Introduction}

Many series exist for $\pi$ and $\pi^2$ \cite{Borwein1987,Pain2022a,Pain2022b}. However, there are only very few published expressions for $\pi^3$. Besides the understanding of such a scarcity from the mathematical point of view and the derivation of summations in order to point out connections between mathematical constants, $\pi^3$ is encountered in physics for instance in the expression of the equivalent of the effective area (used for electromagnetic antennas) for gravitational-wave antennas, which is a measure of the antenna's ability to gather energy from the incident wave \cite{Lewis1995,Pickover2005}. In Ref. \cite{Sun2015}, Sun mentions that the only well-known series for $\pi^3$ is the following one 
\begin{equation}
\pi^3=32\sum_{k=1}^{\infty}\frac{(-1)^k}{(2k+1)^3}
\end{equation}
and in 2010 \cite{Sun2010}, the same author suggested that
\begin{equation}
\sum_{k=0}^{\infty}\frac{\binom{2k}{k}}{(2k+1)^316^k}=\frac{7\pi^3}{216},
\end{equation}
which was proven later \cite{Pilehrood2010}. In the latter paper, Pilehrood and Pilehrood obtained the following Ap\'ery-like series:
\begin{equation}
\pi^3=32\sum_{k=0}^{\infty}\frac{\binom{2k}{k}}{16^k(2k+1)^3}-24\sum_{k=0}^{\infty}\frac{\binom{2k}{k}}{16^k(2k+1)}\sum_{m=0}^{2k-1}\frac{1}{(2m+1)^2}.
\end{equation}
In the above mentioned article \cite{Sun2015}, Sun derived the following expression
\begin{equation}
\sum_{k=1}^{\infty}\frac{2^kH_{k-1}^{(2)}}{k\binom{2k}{k}}=\frac{\pi^3}{48},
\end{equation}
where $H_n^{(2)}$ denotes, for $n\in\mathbb{N}^*$, the Harmonic number
\begin{equation}
H_n^{(2)}=\sum_{k=1}^n\frac{1}{k^2}.
\end{equation}
Gupta published new series representations of $\pi$, $\pi^3$ and $\pi^5$ in terms of Euler numbers and $\pi^2$, $\pi^4$ and $\pi^6$ in terms of Bernoulli numbers \cite{Gupta2017}. He found the following relation
\begin{equation}
\pi^3=\sum_{n=1}^{\infty}\frac{(-1)^{n+1}}{(2n-1)^3\left(2^{2k+2}-1\right)}2^ {2k+4}(2k+3)!\sum_{j=0}^k\left[-\frac{1}{(2n-1)^2\pi^2}\right]^j\frac{1}{(2k-2j+1)!}
\end{equation}
where $j,k$ are integers and $\geq 0$. Following his success in discovering a new formula for $\pi$, Simon Plouffe \cite{Plouffe2006,Chamberland2011,Berndt} postulated several identities which relate either $\pi^m$ or $\zeta(m)$ to three infinite series. Letting
\begin{equation}
S_n(r)=\sum_{k=1}^{\infty}\frac{1}{k^n(e^{\pi r k}-1)},
\end{equation}
the first two examples are
\begin{equation}
\pi=72~S_1(1)-96~S_1(2)+24~S_1(4)
\end{equation}
\begin{equation}
\pi^3=720~S_3(1)-900~S_3(2)+180~S_3(4).
\end{equation}
In the present work, we show that using a trigonometric series obtained by Euler and involving the functions $x\longmapsto \cot x$ as well as $x\longmapsto \mathrm{cosec}~ x=1/\sin x$, it is possible to derive series expansions for $\pi^3$. We present two of them which are of particular interest since they involve the golden ratio. 

\section{New series for $\pi^3$ involving the golden ratio}

In the book by Borwein and Borwein \cite{Borwein1987}, the following formula, due to Euler
\begin{equation}\label{euler}
\pi^3\left[\cot(\pi x)~\mathrm{cosec}^2(\pi x)\right]=\sum_{n=-\infty}^{\infty}\frac{1}{(x-n)^3}
\end{equation}
is presented (13.b, p. 382). It can be used to derive series expansions for $\pi^3$. In particular, for $x=1/5$ or $x=1/10$ for instance, it is possible to obtain expressions of $\pi^3$ as series multiplied by a coefficient involving the golden ratio. Indeed, let us consider first the case $x=1/5$. One has
\begin{equation}
\cos\left(\frac{\pi}{5}\right)=\frac{\phi}{2}
\end{equation}
and
\begin{equation}
\sin\left(\frac{\pi}{5}\right)=\sqrt{\frac{5-\sqrt{5}}{8}},
\end{equation}
where
\begin{equation}
\phi=\frac{1+\sqrt{5}}{2}
\end{equation}
is the golden ratio. One thus has
\begin{equation}
\cot\left(\frac{\pi}{5}\right)=\frac{\sqrt{2}\phi}{\sqrt{5-\sqrt{5}}}=\frac{\phi}{\sqrt{3-\phi}}
\end{equation}
and
\begin{equation}
\mathrm{cosec}\left(\frac{\pi}{5}\right)=\frac{2\sqrt{2}}{\sqrt{5-\sqrt{5}}}=\frac{2}{\sqrt{3-\phi}}.
\end{equation}
Therefore, using Eq. (\ref{euler}) for $x=1/5$ yields
\begin{empheq}[box=\fbox]{align}
\pi^3=\frac{125}{4}\frac{(3-\phi)^{3/2}}{\phi}\sum_{n=-\infty}^{\infty}\frac{1}{(1-5n)^3}.
\end{empheq}
In the same way, setting $x=1/10$, a similar expansion can be deduced. One has indeed
\begin{equation}
\cos\left(\frac{\pi}{10}\right)=\frac{\sqrt{10+2\sqrt{5}}}{4}=\frac{\sqrt{2+\phi}}{2}
\end{equation}
and
\begin{equation}
\sin\left(\frac{\pi}{10}\right)=\frac{1}{2\phi},
\end{equation}
yielding
\begin{equation}
\cot\left(\frac{\pi}{10}\right)=\phi\sqrt{2+\phi}
\end{equation}
as well as
\begin{equation}
\mathrm{cosec}\left(\frac{\pi}{10}\right)=2\phi
\end{equation}
giving finally the expansion
\begin{empheq}[box=\fbox]{align}
\pi^3=\frac{250}{\phi^3\sqrt{2+\phi}}\sum_{n=-\infty}^{\infty}\frac{1}{(1-10n)^3}.
\end{empheq}
Additional expressions can of course be obtained\footnote{For instance, using $x=1/4$, since $\cot(\pi/4)=1$ and $\mathrm{cosec}(\pi/4)=\sqrt{2}$, one gets $\pi^3=32\sum_{n=-\infty}^{\infty}\frac{1}{(1-4n)^3}$.} for other values of $x$, but only a few of them will only involve the golden ratio. As an example, using $x=\pi/15$, one has
\begin{equation}
\cos\left(\frac{\pi}{15}\right)=\frac{1}{8}\left(\sqrt{30+6\sqrt{5}}+\sqrt{5}-1\right)
\end{equation} 
and 
\begin{equation}
\sin\left(\frac{\pi}{15}\right)=\frac{1}{16}\left(2\sqrt{3}-2\sqrt{15}+\sqrt{40+8\sqrt{5}}\right)
\end{equation} 
which unfortunately involves $\sqrt{3}$...

\section{Conclusion}

We proposed two representations of $\pi^3$ involving infinite series and the golden ratio. Although the number of expressions of this type is probably rather limited, the technique can be easily applied to derive other series, relating $\pi^3$ to other mathematical constants. It is worth mentioning that the third power of $\pi$ finds applications in several fields of physics, such as gravitational-wave antennas, in the expression of the effective area.


\begin{thebibliography}{99}

\bibitem{Borwein1987} J. M. Borwein and P. B. Borwein, {\it Pi and the AGM}, John Wiley and Sons, New York, 1987.

\bibitem{Pain2022a} J.-C. Pain, {\it A double series for $\pi$ using Fourier series and the Grothendieck-Krivine constant} (2022),\\ 
\url{https://arxiv.org/abs/2206.05610}

\bibitem{Pain2022b} J.-C. Pain, {\it Relations between $\pi$ and the golden ratio $\phi$ in the form of Bailey-Borwein-Plouffe-type formulas} (2022),\\ 
\url{https://doi.org/10.48550/arxiv.2205.08617}

\bibitem{Lewis1995} M. A. Lewis, {\it Gravitational-wave versus electromagnetic-wave antennas}, IEEE Antennas and Propagation Magazine {\bf 37}, 26-31 (1995).

\bibitem{Pickover2005} C. A. Pickover, {\it A passion for mathematics} (John Wiley \& Sons, Inc., Hoboken, New Jersey, 2005).

\bibitem{Sun2015} Zhi-Wei Sun, {\it A new series for $\pi^3$ and related crongruences}, Int. J. Math. {\bf 26}, 1550055 (2015).

\bibitem{Sun2010} Z. W. Sun, {\it Conjecture on a new series for $\pi^3$}, A Message to Number Theory List (sent on March 31, 2010).

\bibitem{Pilehrood2010} K. H. Pilehrood and T. H. Pilehrood, {\it Series acceleration formulas for beta values}, Discrete Math. Theor. Comput. Sci. {\bf 12}, 223-236 (2010).

\bibitem{Gupta2017} H. C. Gupta, {\it New series representations of $\pi$, $\pi^3$ and $\pi^5$ in terms of Euler numbers and $\pi^2$, $\pi^4$ and $\pi^6$ in terms of Bernoulli numbers} (2017),\\ \url{https://arxiv.org/ftp/arxiv/papers/1710/1710.04083.pdf}

\bibitem{Plouffe2006} S. Plouffe, {\it Identities inspired by Ramanujan notebooks} (part 2), April 2006,\\
\url{http://www.plouffe.fr/simon/}

\bibitem{Chamberland2011} M. Chamberland and P. Lopatto, {\it Formulas for odd Zeta values and powers of $\pi$}, Journal of Integer Sequences {\bf 14}, 236147 (2011).

\bibitem{Berndt} B. Berndt, {\it Ramanujan Notebooks} (volumes I to V), Springer Verlag, 1985 to now.

\end{thebibliography}
\end{document}